\documentstyle[12pt]{article}
 
\newcommand{\R}{{\bf R}}
\newcommand{\E}{{\bf e}}

\newcommand{\PE}{{\bf e_p}}

\newcommand{\Z}{{\bf Z}}
\newcommand{\C}{{\bf C}}

\newcommand{\Om}{{\Omega}}
\newcommand{\om}{{\omega}}
\newcommand{\eps}{{\varepsilon}}
\newcommand{\de}{{\delta}}

\newcommand{\be}{{\beta}}
\newcommand{\g}{{\gamma}}
\newcommand{\Ga}{{\Gamma}}
\newcommand{\la}{{\lambda}}
\newcommand{\La}{{\Lambda}}
\newcommand{\ka}{{\kappa}}
\newcommand{\si}{{\sigma}}
\newcommand{\Si}{{\Sigma}}
\newcommand{\Ham}{{\rm Ham }}
\newcommand{\im}{{\rm Im }}
\newcommand{\Mm}{{\cal M }}
\newcommand{\Nn}{{\cal N }}

\newcommand{\Jj}{{\cal J}}

\newcommand{\proof}[1]{\noindent{\bf Proof#1:\  }}
\newcommand{\jdef}[1]{{\bf #1}}

\newcommand{\MS}{{\medskip}}
\newcommand{\NI}{{\noindent}}
\newcommand{\QED}{\hfill$\Box$\medskip}
\newtheorem{theorem}{Theorem}[section]
\newtheorem{cor}[theorem]{Corollary}
\newtheorem{remark}[theorem]{Remark}
\newtheorem{lemma}[theorem]{Lemma}

\newtheorem{prop}[theorem]{Proposition}
 
\begin{document}

\title{The Geometry of Symplectic Energy} \author {Fran\c{c}ois Lalonde
\thanks{Partially supported by grants NSERC OGP 0092913
 and FCAR EQ 3518} \\UQAM, Montr\'{e}al  \and Dusa McDuff\thanks{Partially supported by
NSF grant DMS 9103033} \\ SUNY, Stony Brook}
\maketitle
 
\section*{Introduction}
 
One of the most striking early results in symplectic topology is 
Gromov's
\lq\lq Non-Squeezing Theorem''
which says that it is impossible to embed a large ball
symplectically into a thin cylinder of the form $\R^{2n} \times B^2$, where
$B^2$ is a $2$-disc.
  This led to
Hofer's discovery of symplectic capacities, which give a way of measuring the size of
subsets in symplectic manifolds.  Recently, Hofer found a way to measure the size (or
energy) of symplectic diffeomorphisms by looking at the total variation of their
generating Hamiltonians.   This gives rise to a bi-invariant (pseudo-)norm on the group
$\Ham(M)$ of compactly supported  Hamiltonian symplectomorphisms of the manifold $M$.  
The deep fact is that this pseudo-norm is a norm; in other words, 
the only symplectomorphism
on $M$ with zero energy is the identity map.
Up to now, this had been  proved only for sufficiently nice symplectic
manifolds, and by rather complicated  analytic arguments.

In this paper we consider a more geometric version of this energy, which was first
considered by Eliashberg and Hofer in connection with their study of the extent to
which the interior of a region in a symplectic manifold determines its boundary.  We
prove, by a simple geometric argument, that both versions of
energy give rise to  
genuine norms  on all symplectic manifolds. Roughly
speaking, we show that if there were a symplectomorphism  of $M$ which
had \lq\lq too little" energy, one could embed a large ball into a thin cylinder
 $M \times B^2$.  Thus there is a direct geometric
relation between symplectic rigidity and energy.
 
The second half of the paper is devoted to a proof of
the Non-Squeezing theorem for an arbitrary manifold $M$.  We do not need 
to restrict to manifolds in which the theory
of   pseudo-holomorphic curves behaves well.    This is of interest
since most other deep results in symplectic topology are
generalised from Euclidean space to other manifolds by using this theory, and hence are
still not known to be valid for arbitrary symplectic manifolds.
 
\section{The Main Results}
 In \cite{H1}, Hofer defined the energy $\|\phi\|_H$ of a compactly supported
Hamiltonian diffeomorphism  $\phi: (M, \om) \to (M,\om)$ as follows:
$$
 \|\phi\|_H =
{\inf}_H (\sup_{x,t}\,H(x,t) - \inf_{x,t}\,H(x,t)),
$$
where  $(x,t) \in M \times
[0,1]$ and $H$ ranges over the set of all  compactly supported Hamiltonian functions
$H: M \times [0,1] \to \R$ whose symplectic gradient vector fields generate  a time $1$
map equal to  $\phi$.    It is easy to check that, for all $\phi, \psi$,
\MS
  
$\bullet\quad$  $\|\phi\|_H = \|\phi^{-1}\|_H$; 

$\bullet\quad$  $\|\phi\circ\psi\|_H \leq \|\phi\|_H + \|\psi\|_H$; and
 
$\bullet\quad$  $\|\psi^{-1}\circ\phi\circ\psi\|_H = \|\phi\|_H$.
\MS
 
Thus $\|\cdot \|_H$ is a symmetric and conjugation-invariant semi-norm on the group
$\Ham(M)$ of all compactly supported Hamiltonian diffeomorphisms of  $M$, and it
follows that the associated function $\rho_H$ given by:
$$
\rho_H(\phi,\psi) =
\|\phi^{-1}\psi\|_H,
$$
 is a bi-invariant pseudo-metric.  However, it is harder to show
that $\|\cdot \|_H$ is a norm, or, equivalently, that $\rho_H$ is a metric. 
Hofer established this
when  $M$ is the standard Euclidean space, using quite complicated analytical
arguments.   This norm is still rather little understood.  A good introduction to its
properties may be found in \cite{H2,H3}.

In this paper, we will  consider the generalized Hofer semi-norm  $\| \cdot \|$ which is
 defined
as follows. Let $M$ be a symplectic manifold of dimension $2n$. If $\partial M \neq
\emptyset$, define Ham($M$) as the group of all compactly supported Hamiltonian
diffeomorphisms which are the identity near the boundary.
 Consider embeddings $\Phi$ of  the strip $M \times [0,1]$ in the product manifold  $(M
\times [0,1]\times \R , \om + dt\wedge dz)$ which are trivial, i.e. equal to  $(x, t)
\mapsto (x, t, 0 )$, for  $t $ near  $0$ and $1$ and for $x$ outside some compact
subset of ${\rm Int}\, M$,  and are such that all leaves of the characteristic
foliation on the hypersurface  $Q =\im\, \Phi$ beginning on   $M\times \{(0,0)\}$ go
through the hypersurface and reach $M\times \{(1,0)\}$.  The induced diffeomorphism
$\phi$ from  $M = M\times \{(0,0)\}$ to $M= M\times \{(1,0)\}$ is  the \jdef{monodromy}
of  $Q$.   Further, the \jdef{energy} of $Q$ is defined to be the minimum length of an
interval  $I$ such that $Q$ is a subset of the product $M \times [0,1]\times I$.

We
define $\|\phi\|$ to be the  infimum of the  energy of all hypersurfaces $Q$ with
monodromy $\phi$. Since $\phi^{-1}$ is the monodromy of the hypersurface $Q$ when read
in the opposite direction, $\|\cdot \|$ is symmetric.  Further, because  the time $1$ map of
the isotopy generated by the function $H(x, t)$ is exactly equal to the monodromy of
the embedding:
$$
(x, t) \mapsto (x,  t, -H(x,t)),
$$
we find that
$$
 \|\phi\| \leq
\|\phi\|_H
$$
 for all $\phi\in \Ham(M)$.  This semi-norm was first considered in
\cite{EH1}.  It is relevant, for instance, when one is trying to understand the
extent to
which the boundary of a region is determined by its interior, since the boundary can
always be $C^0$-approximated by a sequence of hypersurfaces lying inside the region:
see \cite{EH2}.  Note that the two norms defined here might coincide, since no example is
yet known where they differ.

As in Hofer's proof of the non-degeneracy of $\| \cdot \|_H$ in $\R^{2n}$, we will prove that
$\|\cdot\|$ is a norm on $\Ham (M)$ by establishing an energy-capacity inequality which gives
a lower bound for the disjunction energy of a subset in terms of its capacity. Since
all our arguments will rely on properties of embedded balls, the appropriate capacity
to use in the present context is Gromov's radius  $c$.  Thus  for any subset $A \subset
M$, we define
$$
c(A) = \sup \{u: \mbox{there is a symplectic embedding } B^{2n}(u)
\hookrightarrow {\rm Int}\,A \}.
$$
Here,  we use the notation $B^{2n}(u)$ to denote
the standard ball in standard Euclidean space  $(\R^{2n}, \om_0)$ of capacity $u$ and
radius  $\sqrt{u/\pi}$.   Thus the capacity of a ball of radius $r$ is  $\pi r^2$.
In order to distinguish the standard
 balls in $\R^{2n}$  from their images in $M$, we will reserve the dimensional
upperscript to the former only.
 
The \jdef{disjunction} (or \jdef{displacement}) \jdef{ energy} of  $A\subset {\rm
Int}\,M$ is defined to be:  $$\E(A) =
\inf \{ \|\phi\|: \phi\in \Ham(M), \,\phi(A) \cap A = \emptyset\}. $$
We will also need to consider maps $\phi$ which not only disjoin $A$, but also move
$A$ to a new position which is sufficently separated from the old one.  This gives us the
notion of \jdef{proper disjunction energy}.  This is easiest to define for balls. 
  A disjunction $\phi$ of $B(c)$ is said to be \jdef{proper}
if (some parametrization of) $B(c)$ extends to a ball $B(2c)$ such that 
$\phi(B(c))\cap B(2c) = \emptyset$, and the
\jdef{proper disjunction energy $\PE(B(c))$} is the infimum of the energies of all proper
disjunctions of $B(c)$.  Similarly, $\phi$ is said to be a proper disjunction of $A$ if
each ball $B(c) \subset A$ may be extended to a ball $B(2c)$ such that
$$
\phi(A) \cap (A \cup B(2c)) = \emptyset,
$$
and the proper disjunction energy $\PE(A)$ is the minimum energy of such a $\phi$.
 As usual, if there are no (proper) disjunctions of $A$ in $M$, we define its (proper)
disjunction energy to be infinite.

Our  main result is:
 
\begin{theorem}\label{main}  Let $(M, \om)$ be any symplectic manifold, and $A$
any
compact subset of ${\rm Int}\,M$.  Then

\MS
\NI
(i)$\qquad\quad \PE(A) \geq c(A),\quad$ and

\smallskip
\NI
(ii)$\qquad\quad \E(A) \geq \frac 12 c(A). $
\end{theorem}
 
\begin{cor} For any symplectic manifold $M$, $\|\cdot\|$ is a (non-degenerate) norm on
Ham($M$). Hence $\|\cdot \|_H$ is also non-degenerate. \end{cor}

\begin{remark}\rm {\bf (i)} It is very easy to see that the disjunction energy of a
ball  in $\R^{2n}$ is exactly equal to its capacity.  Indeed, when $n = 1$,  an open
ball can be identified with a square and then disjoined  by a translation of energy
equal to its capacity.  In higher dimensions, the result follows from this by
considering the ball as a subset of a product of squares.  \QED
 
\noindent
 {\bf (ii)} It is also easy to check that any disjunction of a ball in $\R^{2n}$ is a
proper disjunction.  (See the proof of Proposition~\ref{basic}.) However,  very little
is known  about the space $Emb(B^{2n}(u), M)$ of balls of  capacity $u$ in an arbitrary
symplectic manifold $M$.   For example, if $n > 2$  it is not even known whether
$Emb(B^{2n}(u), B^{2n}(u'))$ is path-connected when $u < u'$.   Therefore, even if one
restricts to balls $B$  of capacity less than $c(M)/2$, it is not clear what the
relation is between $\E(B)$ and $\PE(B)$. \QED

\noindent {\bf (iii)}  Our arguments actually prove more than what is stated above,
because they are \jdef{local}: they use only the part of the hypersurface $Q$ swept out
by the characteristics emanating from $A$.  We will say that a piece $Q$ of
hypersurface in $M \times [0,1]\times \R$ disjoins $A\subset M$ if, for all $x \in A$,
the point $(x,0,0)$ is one end of a characteristic on $Q$, the other end of which is  a
point  in $( M-A)\times \{(1,0)\}$.  Then we can define the
energy $\E(Q)$ of $Q$ to be the minimal length of an interval $I$ such that  $Q
\subset
M\times [0,1]\times I$, and our result reads:

\MS
 $\qquad\qquad \E(Q) \geq
\frac 12 {c(A)}\quad$ if $Q$ disjoins $A$, and
\smallskip
 
 $\qquad\qquad\E(Q) \geq c(B)\quad$ if $Q$
properly disjoins the ball $B$.

\MS
\NI
In a similar way, we can define and
estimate the energy of a local Hamiltonian diffeomorphism of $M$.\QED
 
\noindent {\bf (iv)}  Polterovich,
using geometric arguments which are very similar in spirit to ours, established in
\cite{POL} that $\|\cdot\|_H$ is a norm on rational symplectic manifolds which are tame at
infinity.  His result is not as sharp as ours, because
he considered the disjunction energy of Lagrangian 
submanifolds, which are  more
unwieldly than balls.\QED \end{remark}
 
Our methods also allow us to prove Gromov's Non-Squeezing 
Theorem in full generality.

\begin{theorem}[Non-Squeezing Theorem]\label{thm:squ} Let $(M, \om)$ be any symplectic
manifold, and denote by $M\times B^2(\la)$ the product of $M$  with the disc $B^2(\la)$
of area $\la$, equipped with the product form.  Then $$ c(M\times B^2(\la)) \leq  \la.
$$ \end{theorem}
 
\noindent \begin{remark}\rm This result was first proved for manifolds such as the
standard $\R^{2n}$ and $T^{2n}$ by Gromov in \cite{GR}. Its range of validity was extended 
by improvements in the
understanding of the behavior of pseudo-holomorphic curves.  However, this method has
definite limitations and is not yet known to  apply to all manifolds.   (The best result
which can be obtained in this way is described in \S3.)  We manage to overcome these
limitations by using
  the techniques which we developed to prove Theorem~1.1.
 As we shall see in Remark~\ref{12} below,
these two theorems are very closely connected, and we will, in fact,  deduce
Theorem~1.1 from Theorem~1.4. To the authors' knowledge, these are the first deep
results in symplectic topology which have been established for all symplectic manifolds.
Arnold's conjecture, for example, has still not been proved, even for all compact
manifolds. \end{remark}

Finally, we observe that the methods developed here permit the construction of some
new embeddings of ellipsoids into balls.  In particular, it is possible to
 solve a problem posed by
Floer, Hofer and Wysocki in \cite{FHW}. This is discussed further in Remark~\ref{ellip}. 
 
Throughout this paper, all embeddings and isotopies will always be assumed to preserve
the symplectic forms involved.
 
We are grateful to Yakov Eliashberg and Leonid Polterovich for useful discussions on
some basic ideas developed in this paper and to Lisa Traynor for showing us explicit
full  embeddings of balls that inspired our Lemma~\ref{Lisa}. The first author thanks
Stanford University for a stay during which part of this work was undertaken.

\section{The energy-capacity inequality in $\R^{2n}$.}
\MS
 
This section  presents a very simple proof of the energy-capacity inequality for
 subsets of
$\R^{2n}$.  The basic idea is that if a hypersurface $Q$ of small energy disjoins a
large ball in $M$, one can construct an embedded ball in the product $M\times B^2$
whose capacity is larger than the area of the $B^2$ factor.  But the 
Non-Squeezing Theorem
 states that when $M = \R^{2n}$  this area is an upper bound for the capacity of
embedded balls in $M\times B^2$.

We need an auxilliary lemma about decompositions of the ball, which was inspired by
Traynor's constructions in \cite{TR}.  Given any set $A$, we will write $\Nn(A)$ to
denote some small neighborhood of it.
 
\begin{lemma}\label{Lisa}  Suppose that $0 < c < C$, and let  $Y = P_1\cup L \cup
P_2\subset \R^2$ be the union of two rectangles, $P_1$ of area $C - c$ and $P_2$ of
area $ c$, joined by a line segment  $L$.  Further, let
$$
Z_{C,c}
=B^{2n}(C) \times P_1 \;\,\cup \;\,B^{2n}(c) \times (L \cup P_2).
$$
Then, there is a
symplectic embedding $B^{2n+2}(C) \hookrightarrow \Nn(\,Z_{C,c})$ in any 
neighborhood of $Z_{C,c}$.
\end{lemma}
\proof{}  Let $\pi$ be the projection $B^{2n+2}(C)\to B^2(C)$ which is induced by
projection onto the last two coordinates.  This represents $B^{2n+2}(C)$ as a kind of
fibration over the disc, with fibers which are concentric balls of different capacities.
Note that the set
$$
 \{\,x \in B^2(C): c(\pi^{-1}(x)) \geq c \}
$$
is exactly $ B^2(C- c)$.  It is easy to see that there is an area preserving embedding $g:
B^2(C) \hookrightarrow \Nn(Y)$ which takes $ B^2(C-c)$ into a neighborhood of $P_1$.  
In fact, we may choose $g$ so that it sends an open neighbourhood of
$$
B^2(C-c) \: \cup \:  (\, (-\infty,0) \times \{0\} \, \cap
\,  B^2(C) \, ) \: \subset \: B^2(C)
$$
into a neighborhood of $P_1$.  Clearly, $g$ is covered by the desired embedding of $B^{2n+2}(C) $
into
$\Nn(Z_{C,c})$.  \QED

\begin{prop}\label{basic} For any compact subset $A \subset \R^{2n}$, 
$$ \E(A) \geq c(A). $$
\end{prop}
\proof{}  Let  $Q$ be a hypersurface of energy $e$ which disjoins
$A$, and let $B\subset A$ be the image of a standard ball of capacity $c$.  We must
show that $e \geq c$.  This will follow if, for any $\de > 0$, we can find an embedding
of the ball $B^{2n+2}(2c)$  of capacity $2c$ into the product $\R^{2n} \times B^2(e +
c + \eps)$,  since the Non-Squeezing Theorem then tells us  that
$$
2 c   \le e +   c + \eps.
$$
 
By Lemma~\ref{Lisa}, it  suffices to embed  $Z_{2c,c}$ symplectically into $\R^{2n} \times
X$, where $X$ is an annulus of area $e+c$.   By hypothesis, there is a rectangle $R$
in $[0, 1]\times \R$ of area $e$ such that $Q\subset\R^{2n} \times R$.  Note that, by
hypothesis, $Q$ is flat  near its ends, that is, that  $Q$ coincides with the
hypersurface  $\{z = 0\}$ near the boundary  $t = 0,1$.  (Recall that we use the
coordinates $(t,z)$ on $[0,1]\times \R$.)    Let  $R'$ be another rectangle in $\R ^2$
of area $c$ with one edge along  $t = 1$, chosen  so that $ R \cup R'$ is a rectangle
of area $e + c$ with one edge along  $t = 0$ and another along  $t = t_1 > 0$. Then form
$X$ by identifying these two edges. 

Let $g:B^{2n}(c) \to B \subset A$ be a symplectic embedding, and extend $g$ to an
embedding, which we also call  $g$, of $B^{2n}(2c)$ in  $\R^{2n}$.  This is possible
because the space of embedded symplectic balls of any given radius in $\R^{2n}$ is path
connected.  (The space of embedded balls of variable radius in any manifold $M$ is
always connected, so long as $M$ is connected, and, when $M = \R^{2n}$ we can fix the
size of the radius by composing with appropriate homotheties.  A similar argument
shows that the space of  symplectic embeddings of two balls in $\R^{2n}$
is path-connected.)  This implies that  any ball in $\R^{2n}$ is isotopic to a standard
one and thus can be extended as much as we wish. Further, we may suppose that the ball
$g(B^{2n}(2c))$ is disjoint from $\phi_Q(B)$, where $\phi_Q$ is the monodromy of $Q$.
To see this, note that because the space of symplectic embeddings of two balls in
$\R^{2n}$ is path-connected, there is a symplectomorphism $\tau$ which is the identity
on $B$ and which moves $\phi_Q(B)$ far away.  Hence we may alter $Q$ without
changing its energy to a hypersurface  with the conjugate monodromy
$\tau^{-1}\circ\phi_Q\circ\tau$.
 
We now define the embedding  $Z_{2c,c} \to \R^{2n}\times X$ as follows.
\MS
 
$\bullet$   $B^{2n}(2c) \times P_1$ goes to $\R^{2n}\times R'$ by $g\times i$,
where $i: U_1 \to R'$.
\smallskip
 
$\bullet$   $B^{2n}(c) \times L$  maps  to the hypersurface  $Q\subset \R^{2n}\times
R$ by a map which takes each line $\{x\}\times L$ to the corresponding flow line of the
characteristic flow on $Q$.
\smallskip
 
$\bullet$   $B^{2n}(c) \times P_2$ goes onto $\phi_Q(B)\times R'$ by the map
$(\phi\circ g) \times i$.
\medskip
 
It is  easy to check that this map preserves the symplectic form.  Hence, the
symplectic neighborhood theorem implies that it extends to the required
symplectomorphism  from $\Nn(Z_{2c,c})$ to $\R^{2n}\times B^2(e + c + \eps)$. \QED
 
\begin{remark}\label{12} \rm The above argument clearly proves Part (i) of
Theorem~\ref{main} for manifolds for which the Non-Squeezing theorem holds.  For, suppose
that the monodromy $\phi_Q$ of a hypersurface  $Q$
is a proper disjunction of $A$ of energy $e$.  Then, for every ball $B(c)\subset A$,
$\phi_Q$ disjoins  $B$ from a ball $\tilde B(2c)$ of twice the capacity and one can
construct an embedding  $Z_{2c,c} \hookrightarrow M\times B^2(e + c + \eps)$ as above.

Similarly, one can prove Part (ii) of Theorem~\ref{main} in this case, by using a
slightly different embedding.  Note that the ball $B^{2n+2}(c)$ embeds into the
set
$$
W_{c,c} =B^{2n}(c) \times Y
$$
where now  $P_1$ and $P_2$
both  have area $c/2$.  This gives a ball of capacity $c$ in $M \times B^2(e + c/2
+ \eps)$.  Therefore, if the Non-Squeezing Theorem holds, we must have $$ c \le  e + c/2 +
\eps,  $$ for all $\eps > 0$ and hence $e \geq c/2$, as claimed. \QED \end{remark}
 
\begin{remark}\label{ellip}\rm Let $E(c_1, c_2)$ denote the ellipsoid  
$$
\Sigma_i \pi(x_i^2 +
y_i^2)/c_i  \leq 1,
$$
where  $c_1 \leq c_2$.  Floer, Hofer and Wysocki show in \cite{FHW}
that if $c_1  \geq 1/2$,  $E(c_1, c_2) $ embeds symplectically in $B^4(1) = E(1,1)$  
if and only if there is a symplectic linear embedding from $E(c_1, c_2)$ into $E(1,1)$, i.e.
only if $a_2 \leq 1$. They asked whether this is sharp. In other words, if $c_1 < 1/2$, is
there $ c_2 > 1$ such that $E(c_1, c_2)$ embeds in $B^4(1)$?  Our embedding methods allow
us to answer this question in the affirmative.  

The idea is as follows.  As Traynor points out in \cite{TR},  the ellipsoid $E(c_1, c_2)$
may be considered to be fibered over the disc $B^2(c_2)$, with fibers which are
smaller by a factor of $c_1/c_2$ than those of the corresponding ball $B^4(c_2)$.  If
$c_1/c_2 < 1/2$, we can therefore fit two of these fibers in the corresponding fiber of
the  ball.  From this, Traynor constructs a full filling of the ball $B^4(1)$ by two open
ellipsoids. 

Now, suppose that we split $E(c_1, c_2)$ into two, considering it to be
contained in a neighborhood of a set $Z$ such as $ Z_{C,c}$ above.  Then, it is not hard
to construct the desired embedding   $E(c_1, c_2)\to B^4(1)$ by folding the two parts of
$Z$ on top of each other.  Details of this construction will be published elsewhere. \QED 
\end{remark}

\section{The Non-Squeezing Theorem}
 
In this section we will  use the theory of $J$-holomorphic curves to prove  the
Non-Squeezing Theorem under certain hypotheses which look  rather artificial.  We will see
in \S4 how to construct families of embeddings which satisfy them.
 
 A  symplectic manifold $(M,\om)$ is often said to be \jdef{rational}  if the
homomorphism induced by $[\om]$ from $\pi_2(M)$ to $ \R$ has discrete image.  In this
case and if this image is not $\{0\}$, we will call the positive generator of this image
the \jdef{index of rationality} of $(M,\om)$, denoted $r(M)$ or simply $r$. If  the
image is $\{0\}$, we set $r(M) = \infty$ and if $M$ is not rational, this index is set
equal to $0$ (this does not quite follow the conventional definitions).
 
 We will consider $M\times S^2(\la)$ with a product form $\Om = \om \oplus \la \si$
where $\si$ is normalised to have total area $1$, and will say that a ball in $M \times
S^2(\la)$ is \jdef{standard} if it is the image of an embedding of the form $$
B^{2n+2}(c)\hookrightarrow B^{2n}(c) \times B^2(c) \hookrightarrow M \times S^2(\la) $$
where the first map is the obvious inclusion and the second is a product.  Clearly the
capacity of any standard ball in $M\times S^2(\la)$ is bounded above by $\la$.
The main
result of this section says that this remains true for any ball which is isotopic to a
standard ball through large balls.
 
 \begin{prop}\label{squ1} Let  $M$ be closed and have index of rationality $r > 0$, and
suppose that $g_t:B^{2n+2}(c_t) \to M\times S^2(\la)$ is a family of symplectic
embeddings such that $g_0$ is standard and   $$ c_t \geq \sup \{r,\la - r\},
$$ for all  $t$.
Suppose also that the image of $g_t$ misses one fiber $M \times pt$  for all $t$.
Then  $$ c_t < \la $$ for all $t$. \end{prop}
 
We will begin by explaining the  usual proof of Gromov's 
Non-Squeezing theorem via
pseudo-holomorphic curves.  We assume that the reader is familiar with the basics of
the theory of pseudo-holomorphic curves as explained in \cite{GR,ELL} for example. The
most general argument works when the manifold $(V, \Om) = (M\times S^2(\la), \om +
\la\si)$ is semi-positive (or semi-monotone): see \cite{CTB}.  This condition says that, for
generic tame $J$ there are no $J$-holomorphic spheres  with negative Chern number.
It is satisfied by all manifolds of dimension $\leq 6$ and by manifolds for which there is
a constant $\mu \geq 0$ such that
$$
c_1(\alpha) = \mu [\om](\alpha),\quad\mbox{for all } \alpha \in \pi_2(M),
$$
where $c_1$ is the first Chern class of $(M, \om)$.
 
\begin{lemma} Let $g: B^{2n+2}(c) \to  M\times {\rm Int}\,B^2(\la)$ be a symplectic
embedding, and suppose that $(V, \Om) = (M\times S^2(\la), \om + \la\si)$ is
semi-positive.  Then
$$
c < \la.
$$
\end{lemma}
 
\proof{}  Clearly, we may consider $g$ as an embedding into $V$ which misses one
fiber.  Let $J_{B}$ be an almost complex structure tamed by $\Om$  which extends the
image by $g$ of the standard complex structure on $\R^{2n+2}$. In order to show that
$c \le \la$, it is enough, by Gromov's monotonicity argument, to show that there exists a
$J_{B}$-rational curve $C$ of symplectic area smaller or equal to $\la$ passing through
the  center $g(0)=p_{0}=(q_{0},z_{0})$ of the ball $g(B^{2n+2}(c))$.  The reason
is that the part of this curve in $g(B^{2n+2}(c))$ pulls back to a holomorphic curve
$S$ through the center of $B^{2n+2}(c)$.  Since $S$ is minimal with respect to the
usual metric on $B^{2n+2}(c)$, the monotonicity theorem implies that its area is at
least $c$. Thus $$ c \leq \mbox{area}\,S = \int_S\om_0 < \int_C\Om \leq\la,  $$ as
required.
 
To produce the required $J_B$-curve, we argue as follows.  Let $\Jj $ denote the
 space
 of all $C^\infty$-smooth almost complex structures tamed by  $\Om$ and, for
 $A \in H_{2}(V; \Z)$, let
 $\Mm(A)$  be the space of all pairs $(f,J)$ such that
 $f:\C P^{1} \longrightarrow V $ is $J$-holomorphic,
 $[f] \in A$, and $f$ does not factor through a self-covering.  A basic fact in this
respect is that the projection map  $P_A: \Mm(A) \to \Jj$ is Fredholm with index
$2(c_1(A) + n + 1)$, where $c_1$ is the first Chern class of the tangent bundle of $(V,
J)$.   (For more details of this step see \cite{ELL,MS}.)
 
Take any almost complex structure $J_{1}$ tamed by $\om$ on $M$, integrable  in a
neighbourhood $U$ of $q_{0}$,  take the usual complex structure $J_0$ on $S^2$, and
denote by $J_{spl}$ the split structure $J_{1} \oplus J_{0}$.   It is easy to see that
for all $q \in U$ the rational curves $\{q\} \times S^{2}$ in class $A_0=[\{pt\} \times
S^{2}]$ are regular for this complex structure  in the Fredholm sense.   (This means
that, given any holomorphic parametrization  $f$ of these curves, the points $(f,
J_{spl})$ are regular points of $P_{A_0}$.) Note also that $P_{A_{0}}^{-1}(J_{spl}) =
\{q \times S^{2}:q \in M\}$, because any holomorphic map from $\C P^{1}$ to the split $M
\times S^{2}$ in class $A_{0}$  induces, by projection on
 the first factor, a holomorphic map to $M$, which is null-homologous and hence
constant.        Obviously all points  $p \in V$ that project to $U$ are regular
 values of
the evaluation map \[  ev: P_{A_{0}}^{-1}(J_{spl}) \times _{G} \C P^{1} \longrightarrow
V, \] where $G$ is the conformal group of $\C P^{1}$, and these $p$ have pre-image
$ev^{-1}(p)$ containing exactly one point.  It follows that  there exists
 a structure $J'$ near $J_{spl}$ which is  generic, that is  regular for all
projections $P_{A}$, and  is such that some point $p' \in U \times S^{2}$  is still a
regular value for the evaluation map on $P_{A_{0}}^{-1}(J') \times _{G} \C P^{1}$, with
exactly one point in its pre-image.
 
Now, let $J''$ be a generic almost complex structure in the neighbourhood of $J_{B}$
and $\Ga$ a path  in $\Jj $ from $\Ga (0) = J'$ to $\Ga (1) = J''$, transverse to all
projections $P_{A}$. Then $P_{A_{0}}^{-1} (\Ga )$ is a smooth manifold, and we consider
a short path $ \g $ in $V$ from
 $\g (0) = p'$ to a point $\g (1) = p''$ in a neighbourhood of $p_{0}$, such that the
obvious evaluation map
 $$
 ev_{\Ga,A}: P_{A}^{-1}(\Ga ) \times _{G} \C P^{1} \longrightarrow V \times [0,1]
$$ 
is
transverse to $\g \times id :[0,1] \rightarrow V \times [0,1]$ for all $A$. (Here 
$ev_{\Ga,A}$ maps onto $[0,1]$ by projection through $\Ga$.) Denote by $N$ the
one-dimensional submanifold
 $ev_{\Ga,A_{0}}^{-1}(\g \times id)$.  By construction, there is exactly one point of
$N$ which maps to $(\g(0), 0)$. Therefore, if $N$ were compact, it would have at least
one other point over $(\g(1), 1)$.  In other words, there would exist a $J''$-rational
curve in class $A_{0}$ passing through $p''$.  And if this were also true for a
sequence of paths in $\Jj$ whose end points converge to $J_{B}$ and a  sequence of
paths in $V$ whose end points converge to $p_{0}$, this would give, by   Gromov's
compactness theorem (see \cite{GR}), a sequence of holomorphic spheres (weakly)
converging to a $J_{B}$-cusp-curve.  The component of this cusp-curve passing through
$p_{0}$ would have area smaller or equal to $\la$ and so would be the desired
$J_B$-curve.

 Let us suppose now that one of these manifolds $N$ is not compact.  There would
 then
exist a sequence of $J_{i}$-curves $f_{i}:\C P^{1} \rightarrow V$ with $\{J_{i}\}$
converging to some  $\overline{J} = \Ga (t_{0})$ and  $\{f_{i}\}$ diverging. Since $V$
is compact, the compactness theorem implies that some subsequence would converge weakly
to a cusp-curve passing through $\overline{p} = \g (t_{0})$. This cusp-curve would be a
connected union of  $\overline{J}$-curves in classes $A_1, \dots, A_k$, where $$ A_0 =
A_1 + \dots + A_k. $$ Therefore, the proof may be finished if we put a hypothesis on
$M$ which ensures that a generic path $(\Ga,\g)$ does not meet any such cusp-curve.
For example, since $\Om(A_j) > 0$ for $j = 1,\dots, k$, it is clearly enough to assume
that $\pi_2(M) = 0$ or more generally that $\la = \Om(A_0) \leq r$.   The real trouble
comes from the possible presence of multiply-covered curves of negative Chern number,
and it is shown in \cite[\S 4]{CTB} that it suffices to assume that $V$ is semi-positive.
\QED
 
{\bf Proof of Proposition~\ref{squ1}}
\MS
 
There is no semi-positivity hypothesis here: we get around the problem caused 
by cusp-curves by considering  a very special path $\Ga$.  First, let us
consider a path  from $J_0 = J_{spl}$ to $J_1$ such that, at each time 
$t$, $J_t$ is equal to the push-forward by the embedding $g_t$ of the standard
structure on $\R^{2n+2}$.  By assumption on the embeddings,
we may also suppose that one fiber
$M\times pt$ is $J_t$-holomorphic, for each $t$.
 Suppose that there were a $J_t$-holomorphic
$A$-cusp-curve through the center $g_t(0)$ of the ball  for some $t$ with
homology decomposition
$$
A = A_1 + \dots + A_k.
$$
Let
$C$ be the component of this cusp-curve through $g_t(0)$.  
We may suppose that $[C] = A_1.$  The argument
at the beginning of Lemma 3.2 shows that $\Om(A_1) = 
\int_C \Om > c_t$.  Hence, 
$$
\Om(A_j) < \la - c_t \leq r,\quad \mbox{if}\quad j> 1.
$$
Thus the classes $A_j, j>1$, do not lie in $H_2(M)$.  On the other hand,
by positivity of intersections (see \cite{GR,LB}), the fact
that a fiber is $J_t$-holomorphic implies that the
intersection number $A_j \cdot [M]$ is $ \geq 0$
for all $j$.  It follows that one of the $A_j$ has the form $A - B_j$ for
some $B_j \in H_2(M)$, and that the others are all elements of $H_2(M)$.
Putting all this together, we see that the decomposition must have the form
$$
A = B + (A-B),
$$
for some $B \in H_2(M)$ such that $\Om(B) > c_t$.
 
 Note that neither component of this cusp-curve can be multiply-
covered.  For if the component in class $B$ were a $k$-fold covering for
some $k \geq 2$,
we would have to have $\om(B) > 2c_t$, which is impossible 
by assumption on $c_t$ since we must
have $\om(B) < \la$.  Further, the component in class $A-B$ cannot be
multiply covered since $A - B$ is not a multiple class.  

This argument shows that, for all the elements $J_t$ in our 
special path, the only $J_t$-holomorphic
$A$-cusp-curves are of type $(B,A-B)$.  By
the compactness theorem, a similar statement must hold for every
$J$ in some neighborhood of this path.  Thus we may assume that the
points on the regular path $\Ga$ considered above have this property.
The arguments in \cite{CTB} show that these cusp-curves are well-behaved,
and fill out a subset of $V$ of codimension at least $2$.  Therefore,
a generic path $(\Ga,\g)$ will not meet these cusp-curves,
and the argument may be finished as before. \QED
 
The next lemma is the key step in extending our results to non-compact manifolds.
 
\begin{lemma}~\label{noncompact}  Let $M$ be a non-compact symplectic manifold of
index of rationality $r>0$. For each compact subset $K$ of $M$, there is a number $\zeta
> 0$ such that whenever  $g_t:B^{2n+2}(c_t) \to K\times S^2(\la)$ is a family of
symplectic embeddings  missing one fiber $M\times pt$  for all $t$ and beginning with a
standard embedding $g_0$, then
$$
 c_t \geq \la - \min(r, \zeta) \quad\mbox{for all  }
t \quad  \Rightarrow\quad c_t < \la \quad\mbox{for all  } t.
$$
\end{lemma} \proof{}
Let $K_1$ be a compact subset of $M$ which contains $K$ in its interior.  In order to
make the previous argument go through, we just have to ensure that 
the $A$-curves in $N$ do
not escape outside $K_1 \times S^2$. Suppose that $\partial K_1 = \Si$ is a smooth
hypersurface and let $U$ be a compact neighbourhood of $\Si$ disjoint from $K$.
Since
all balls lie in $K$ we may assume that the special almost
complex structures $J_t$ are all equal to the $\Om$-compatible
split structure $J_{spl}$ on $U$.  Because $U$
is compact, it is easy to see that there exists $\zeta, s$ both small enough so that,
given any $J_t$\/-holomorphic curve $C$ 
passing through
some point $p \in \Si$,
the  ${\Om}$\/-area of $C \cap B_s(p)$ is larger than $3\zeta$. Choosing now, in
the proof of Proposition~\ref{squ1}, the generic path of almost complex structures so
that each be sufficiently close to $J_{spl}$ on $U$, we get a lower bound equal to
$2\zeta$ on the area of $C \cap B_s(p)$.  
If $C$ is a curve in the path $N$, 
 the $\Om$-area of the part of $C$ which lies in the ball $\im\,g_t$ must be
at least $c_t$. But $c_t + 2\zeta > \la$, by hypothesis.  Therefore, none
of the $A$-curves in $N$ meet $\Si$. \QED

\section{Embedding balls along monodromies}

By Remark~\ref{12}, all theorems will be proved if we show that the
Non-Squeezing  Theorem holds
for all manifolds. Our tool to do this is Proposition~\ref{squ1}.
Thus, given
a ball  $B(c)$ in a cylinder $M\times B^2(\la)$ with $c > \la$,
we aim to construct
another  ball $B''(c)$ of capacity $c$, which is contained  in some cylinder
$M'' \times B^2(\la'')$  with $c >
\la''$ and which is  isotopic to a standard ball through
large balls.  Proposition~\ref{squ1}
then implies that the ball $B''(c)$ cannot exist, and it follows that $B(c)$
does not exist either.
 
 It is quite a delicate
matter to obtain a ball $B''$ with the required properties, and as our
notation implies we do this by a two-step process, first constructing an intermediate ball
$B'$, and then using that to get $B''$.  We begin by explaining the basic
procedure which constructs these balls, and then will give the proof.
 
\subsection{The $N$\/-fold wrapping construction}
 
\noindent
Because  the original ball $B$ may not extend to a ball of capacity $2c$,
we must use a multiple wrapping process to maintain the capacity of our
balls. Therefore, instead of using the set  $Y$ of Lemma~\ref{Lisa},
we will use the sets $Y_N\subset \R^2$ described below.
It will be often
convenient to use rectangles rather than discs. As always, the label of a
 set will indicate its capacity (or area), and as before we will distinguish the
standard balls (the domain of our maps) from their images by reserving the dimensional
upperscript to the former only.
 
     Let $V^{2m}$ be any symplectic manifold, $g: B^{2m}(\ka) \to V$
a symplectic embedding of a ball of capacity $\ka$, whose image is
denoted $B(\ka)$, and $\phi_s, 0 \leq s < \infty$, a diffeotopy of $V$
such that $\phi_s$ is periodic in $s$ with period $1$.
We assume that $\phi_s$ has
a $1$\/-periodic generating Hamiltonian $H:V \times [0, \infty) \to \R$ which
satisfies
 
$\bullet$ $H$ vanishes near any integral value of $s$;
 
$\bullet$ for each $s$, $\min_{V} H_s = 0.$
 
When referring in this section to the energy of such a diffeotopy $\phi_s$,
we will always mean the maximum over $s \in [0,1]$ of the total variation of
$H_s$. Further, we will say that $\phi_s$ \jdef{strictly disjoins} a ball
$B\subset V$ if  the balls
$$
B, \phi_1(B),\phi_2(B),\dots,\phi_k(B),\dots
$$
are all disjoint.
\medskip
 
Now choose any positive integer $N$, and,  for $i = 1,\dots, N$ define
$a_i, b_i \in \R$ by
$$
a_i = (i-1)(1 + 1/N),\qquad b_i = a_i + 1/N.
$$
Thus $a_{i+1} = b_{i} + 1$.  Let $Y_N = Y_N(\ka) \subset \R^2$ be the union
of $N$ rectangles
$P_i , 1 \leq i\leq N,$ of area $\ka/N$ with $N-1$ lines $L_i, 1 \leq
i\leq N-1,$ of length $1$, where:
$$
P_i = \{(u,v) : a_i \leq u\leq b_i,\; 0 \leq v \leq \ka\},
$$
$$
L_i = \{(u,v): b_i \leq u \leq a_{i+1},\; v = 0\}.
$$

Observe that any
embedding of $Y_N(\ka)$ extends to some neighborhood $\Nn(Y_N(\ka))$ and hence
induces an embedding of the disc $B^2(\ka)$, since this fits inside any
neighborhood of $Y_N(\ka)$.  Therefore an embedding
$$
G: B^{2m}(\ka)\times Y_N(\ka) \hookrightarrow V\times \R^2,
$$
induces an embedding $ G \circ i$ of the ball $B^{2m+2}(\ka)$ by
$$
B^{2m+2}(\ka)\; \stackrel{i}{\hookrightarrow}    B^{2m}(\ka)\times \Nn(Y_N(\ka))
\;\stackrel{G}{\hookrightarrow}\;   V\times \R^2 $$
where $i$ is the obvious inclusion.
 
\medskip

    Now define $G: B^{2m}(\ka)\times Y_N(\ka) \hookrightarrow V\times \R^2$ by
\begin{eqnarray*}
G(p,u,v) & = & (\phi_{i-1}(g(p)), u,v),\quad\mbox{when}\quad (u,v)\in P_i, \\
& = &(\phi_{(i-1) +u - b_i}(g(p)), u, H(\phi_{(i-1) +u - b_i}(g(p)), u-b_i)),
\quad\mbox{when} (u,v)\in L_i.
\end{eqnarray*}
It is easy to check that $G$ is symplectic and so extends to a symplectic embedding (which we
will also call $G$) of  $ B^{2m}(\ka)\times N( Y_N(\ka))$ into $V \times \R^2$.
\medskip
 
The corresponding ball, which is the image of the composite
$$
B^{2m+2}(\ka) {\hookrightarrow} B^{2m}(\ka)\times \Nn(Y_N(\ka)) 
\stackrel{G}{\hookrightarrow} V\times \R^2
$$
will be called the \jdef{unwrapped ball} associated to $B(\ka) = \im\, g$ and $\phi_s$,
and will be denoted $B(B(\ka), \phi)$.  Note that it exists for any isotopy $\phi_s$, not
only for disjoining isotopies.  For the sake of clarity, we will sometimes write $G_\phi$
for the corresponding embedding $G$.
\medskip
 
If $e$ is the energy of $\phi_s$, then $H$ takes values in $[0, e]$ and the image $\im\, G$
of $B^{2m} (\ka)\times Y_N$ by $G$ lies in the set
$$
V \times \left(\coprod P_i\cup\coprod Q_i\right) = V \times S \;\subset \; V\times \R^2,
$$
where $Q_i, i = 1, \dots, N-1$ are the rectangles
$$
Q_i = \{(u,v): b_i \leq u \leq a_{i+1},\; 0 \leq v  \leq e \}.
$$
 
     Denote by $\tau_S$ the
translation of $V \times S$ by $1 + 1/N$ in the $u$\/-direction of $\R^2$.  (We will
make no distinction between the translation of $S\subset \R^2$ and its lift to the
product $V \times S$.)
This sends $P_i$ to $P_{i+1}$ and  $Q_{i}$ to $Q_{i+1}$, and it is easy
to see that, if $\phi_s$ strictly disjoins $B(\ka)$, then  $\tau_S$  strictly disjoins $\im\, G$,
that is, all the balls
$$
\im\, G, \tau_S(\im\, G), (\tau_S)^2(\im\, G), \dots
$$
are disjoint.  Thus,
if $X$ is the annulus obtained by quotienting $S$ by the translation
$\tau_S$, $X$ has area $A = \ka/N + e $ and, as in Proposition~\ref{basic}, we get an
embedded ball in  $V \times B^2(A + \eps)$ by the composite:
$$
B^{2m+2}(\ka) \hookrightarrow B^{2m}(\ka) \times \Nn(Y_N(\ka))
\stackrel{G}{\hookrightarrow} V \times S \to V\times X \hookrightarrow V \times
B^2(A + \eps). 
$$
This ball wraps $N$ times round the annulus $X$, and will be called the \jdef{wrapped ball}
$B_W(B(\ka),\phi)$  generated by $B(\ka)$ and $\phi_s$.
\medskip
 \begin{remark}\label{rem}\rm\quad Note that if $e < \ka$, we may, by choosing $N$
large enough, arrange that $A$ be arbitrarily close to $e$.
Thus, for sufficently small  $\eps > 0$, we get a ball $B_W = B_W(B, \phi)$ of capacity
$\ka$ inside a cylinder in $V \times \R^2$ of area $A +\eps < e + 2\eps <\ka$.
\end{remark}

The next result is obvious.
 
\begin{lemma}\label{iso}   If $B_t = \im\,(f_t: B^{2m}(\ka_t) \to V)$ and
$\{\phi^t\}$ vary smoothly with respect to
a parameter $t$, the corresponding wrapped ball $B_W(B_t, \phi^t)$ varies
smoothly. \end{lemma}
\medskip

\begin{lemma}\label{extend} The translation $\tau_S$ of $V\times S$ which  disjoins
the unwrapped ball $B(B, \phi)$  may be extended to a $1$\/-periodic
diffeotopy $\{\si_s\}_{0\leq s< \infty}$   of $V\times \R^2$ in such a manner that
 
$\bullet$  $\si_1 = \tau_S$;
 
$\bullet$       the diffeotopy
$\si_s$ strictly disjoins $B(B, \phi)$; and
 
$\bullet$ the energy of $\si_s$ is  $\leq\, A + \eps$.
\end{lemma}
\proof{} Suppose first that the rectangles $P_i$,
$Q_i$ in $T$ all have the same $v$-height.  Then $S$ has the form $\R \times I$,
 and one
can extend $\tau_S$  to have the form $(x,u,v)\mapsto (x, u + \be(v), v)$ on
$V \times \R\times I$, for some suitable
 bump function
$\be$ which equals $1 + 1/N$ on the interval $I$.  This map has the generating
 Hamiltonian
$H(x,u,v) = \int^v \be$
which has energy $\leq \int \be$.  The general result follows because there
is an area-preserving map which commutes with $\tau_S$ and takes $S$ into a set $S_0$
of the form $\R \times I$ with area $S_0/\tau < A + \eps$.
\QED

Of course, the wrapped ball $B_W(B, \phi)$ may also be strictly disjoined by a
diffeotopy (a translation)  of energy $<\,A + \eps <\, \ka$. But for the first step of our
argument we consider instead $B(B, \phi)$ with disjoining isotopy $\si_s$ since the
latter  is more flexible.
 
\subsection{Regularity and the plan of the proof}

    We can now make more precise the plan briefly outlined in the introduction of the
section.  We fix a small constant $\eps > 0$ of size to be determined later.  We start
with the ball
$$
B^{2n+2}(c) \;\stackrel{g}{\hookrightarrow} \;B(c) \;\subset\; M \times B(\la) 
\;\subset\; M\times \R^2\; = \;M'
$$
 (with $\la < c$)
given in the statement of the Non-Squeezing Theorem.  This  is strictly disjoined by a
diffeotopy $\phi_s$ of energy $< \la + \eps$  which translates points in the $\R^2$
direction.     We first construct the unwrapped ball $B'(c) = (B(c),\phi) \subset
M'\times \R^2$, together with the disjoining diffeotopy $\si_s$.    By
Lemma~\ref{extend} above, we may choose $N$ so large that $\si_s$ has energy $<\; \la
+ 2\eps$. Thus the wrapped ball $B''(c) = B_W(B'(c), \si_s)$ of capacity $c$ lies in a
cylinder $C'' = M' \times \R^2\times B^2(\la + 3\eps)$.
 
\begin{prop}\label{B''}  The ball $B''(c)$ in $C'' =  M' \times \R^2\times B^2(\la + 3\eps)$
is isotopic through balls of capacity $\geq \la$ to a standard ball.
\end{prop}
 
\begin{cor} The Non-Squeezing Theorem holds for any symplectic manifold.
\end{cor}

\NI
{\bf Proof of the corollary:}  If  $(M, \om)$ is compact and rational with index of rationality 
$r$, choose $\eps < r/3$.  The Non-Squeezing Theorem then follows immediately from
Propositions ~\ref{B''} and ~\ref{squ1}: simply apply Proposition ~\ref{squ1}
to the closed manifold $M \times T^2(\La) \times T^2(\La)$ which has the same
index of rationality as $M$, where $\La$ is chosen large enough so that
$M \times T^2(\La) \times T^2(\La) \times B^2(\la + 3\eps)$ contains the ball
$B''(c)$.  If $(M, \om)$ is compact but not rational, one
can slightly perturb both the form $\om$ on $M$ and
the ball $B(c) $ in $M\times B^2(\la)$ to get a ball $\tilde B$ 
sitting inside $\tilde M\times B^2(\la)$, where $(\tilde M,\tilde\om) $ is rational.
Although $\tilde B$ may  have capacity $\tilde c$ a little less than $c$, we may clearly
arrange that if $c > \la$ then  $\tilde c > \la$.  Thus the  Non-Squeezing Theorem
holds true for all compact manifolds.
 
      Suppose now that $M$ is non-compact. Note that the initial ball $B(c)$ sits in a
compact region of $M \times B^2(\la)$, so that we may assume both that $M$ is a
compact manifold with boundary and that it has positive index of rationality $r$.
Then $B(c) \subset K \times B^2(\la)$, where $K$ is a compact subset of
$M$. Hence $B''(c)$ lies in $K \times T^2(\La) \times T^2(\La) \times
S^2(\la + 3\eps)$ where $\La$ increases as $\eps$ decreases (and as $N$
increases), because the $N$\/-fold wrapping construction does not move the ball
in the $M$\/-direction. This ball is isotopic to a standard ball through balls
in $K \times T^2(\La) \times T^2(\La) \times S^2(\la + 3\eps)$ of capacity
equal to $\la + 3\eps$ up to a small quantity $3 \eps$. But note that
$r(M) = r(M \times T^2(\La) \times T^2(\La))$ and that the constraining number
$\zeta$ of Lemma~\ref{noncompact} depends only on $K$: it is independant of
the size $\La$ of $T^2(\La)$ as soon as $\La$ is large enough. Therefore,
for $\eps < \frac{\min (r,\zeta)}{3}$, the argument of Lemma~\ref{noncompact}
applies.       \QED

In order to explain our strategy for proving Proposition~\ref{B''} it will be convenient to
introduce the following definition.
 
Given $\de > 0$, we will say that an isotopy $\phi_s$  in some manifold is \jdef{$\de$-regular} 
on a ball $B(\ka)$ if, for each $s \in [\de, 1]$,
the balls
$$
B(\ka), \;\phi_s(B(s\ka)),\; \phi_{s}^2(B(s\ka)), \ldots, \;\phi_{s}^k(B(s\ka)), \ldots
$$
are all disjoint (note, in particular, that $\phi_s$ strictly disjoins $B(\ka)$).
Here $B(s\ka)$ denotes a ball of capacity $s \ka$ which is concentric  with $B(\ka)$.
Notice that it is not so much the isotopy itself which is important but its relation to
the concentric  balls $B(\ka)$.  Further, all
the balls $\phi_{ks}(B(s\ka))$ are assumed to be disjoint from the whole initial
 ball
$B(\ka)$.
 
 Regularity is really
a $1$\/-dimensional notion: any translation of $\R$ at constant speed
which disjoins a given interval
also disjoins subintervals within a time proportional to their lengths.
The basic $2$\/-dimensional example of a regular diffeotopy
is  a translation which disjoins a rectangle in $\R^2$, or its conjugate which disjoins a
disc. (We will often call the latter a translation too.)
 
To be more precise, let   $\phi_s, 0
\leq s < \infty$, be the translation  of the strip $S = \R \times [0,h]$ in the
$u$\/-direction at speed $\nu > 1$. (As above, we use the coordinates $(u,v)$ on $S$.) Note
that this isotopy is  generated by a Hamiltonian which is a linear function of $v$ on the
strip, and vanishes outside some slightly larger strip.  As in Lemma~\ref{extend}
above, its total variation may be taken arbitrarily close to $h\nu$. Then for any small
$\de >0$ and any $\ka < h$, consider an embedding of the standard disk $g: B^{2}
(\ka) \hookrightarrow S$
such that $B(s\ka) = g(B^2(s \ka))$ lies inside $[-s,0] \times [0,h]$ for
all $s \in [\de,1]$. (Such an embedding exists because our choice of constants
$h > \ka$, $\de > 0$ leaves a little extra room.)  It is easy to check that $\phi_s$ is
regular on $\im\, g$.

      With this simple basic example, we can construct higher dimensional examples of
regular disjoining diffeotopies since regularity is stable under products.

\begin{lemma}~\label{regproduct}  Given a $\de$\/-regular pair $(f, \phi)$ in $\R^2$ and 
$g$ any symplectic embedding of a ball in any manifold, form the product 
$h=g \times f: B^{2n}(\ka) \times B^2(\ka) \to 
V = M \times \R^2$. Then the pull-back $\pi^{\ast}(\phi_s)$ of $\phi$ by the 
projection  $\pi:M\times \R^2 \to \R^2$
is also $\de$\/-regular on the composite  $h \circ i$ of $h$ with the standard embedding 
$i: B^{2n+2}(\ka) \hookrightarrow B^{2n}(\ka) \times B^2(\ka)$. 
\end{lemma}

\proof{} This is clear because 
$i$ sends each concentric subball of capacity $\ka'$ onto a ball in $B^{2n}(\ka) \times B^2(\ka)$
whose projection on the second factor is a concentric subdisc of same capacity $\ka'$. \QED
\medskip

The following lemmas form the heart of our argument. They are proved in the following section.

\begin{lemma}~\label{unwrap}  For every $\de, \eps > 0$, there is a
$1$-parameter family $B_t'(\ka_t), \si^t_s$,$0 \leq t\leq 1$, of balls and strict disjoining
isotopies in $M' \times \R^2$ such that:
\MS
 
$\bullet\quad$   $B_0'(\ka_0) = B'(c)$ and $\si^0_s = \si_s$;
 
$\bullet\quad$  the final isotopy $\si^1_s $ is $\de$-regular over
the ball   $B_1'(\ka_1) $; 
 
$\bullet\quad$ the isotopy $\si^t_s$ has energy $< \la + 2\eps$  for all $t$; and
 
$\bullet\quad$ the balls $B_t'$ have capacity $\ka_t \geq
\la$ for all $t$, and $\ka_1 = \la$.
\end{lemma}

 \begin{lemma}~\label{regiso}  Suppose given a ball $B$ of capacity $\la$ in $V$
 with
$\de$-regular strictly disjoining isotopy $\si_s$ of energy $ <\la + 2\eps$, and
 let $B_W$
be the corresponding wrapped ball in the cylinder $C = V \times B^2(\la + 3\eps)$.
 Then, if $\de$ is sufficiently small, $B_W$ is isotopic through balls in $C$ of
capacity $ \la$ to a standard ball.
\end{lemma}

\noindent
{\bf Proof of Proposition~\ref{B''}} \quad By Lemmas~\ref{iso} and~\ref{unwrap}, the wrapped
ball $B''(c)$ is isotopic through balls of capacity $\ka_t \geq\la$ which embed in
$C''$ to the wrapped ball $B_1'' = B_W(B'_1(\ka_1), \si^1)$.  Since the   isotopy $\si^1_s$
is regular,  Lemma~\ref{regiso} shows that the wrapped ball $B_1''$ is isotopic in $C''$ to
a standard ball.\QED

\subsection{The Construction of Isotopies}
 
Our first lemma constructs an isotopy of an unwrapped ball, and the second one for a
wrapped ball. We begin with the proof of the second lemma.
\medskip
 
\noindent
{\bf Proof of Lemma~\ref{regiso}}$\quad$
The wrapped ball $B_W$ is the image of a composite
$$
B^{2m+2}(\la) \hookrightarrow B^{2m}(\la) \times \Nn(Y_N(\la)) \stackrel{G}{\hookrightarrow} V
\times S \to  V\times B^2(\la + 3 \eps),
$$
where the last map is essentially the quotient by the translation $\tau_S$.
The easiest way to describe the isotopy from $B_W$ to a standard embedding is to
construct an isotopy of the map
$$
B^{2m+2}(\la) \hookrightarrow B^{2m}(\la) \times \Nn(Y_N(\la)) \stackrel{G}{\hookrightarrow}V
\times S $$
whose image is disjoined by all iterates of  $\tau_S$.
 
The
crucial observation is that, just as in Lemma~\ref{Lisa}, the ball $B^{2m+2}(\la)$ fits
inside a neighborhood of the  subset $W$ of $B^{2m}(\la)\times
Y_N(\la)$ which is defined as follows:
\begin{eqnarray*}
 W = \left(B^{2m}(\la)\times P_1\right) \;\;\cup\;\;\left(
B^{2m}(\la(1-1/N)) \times (L_1 \cup P_2)\right) \;\;\cup\;\;\dots  \vspace{7mm}\\
\hspace{10mm} \cup \;\;\left(B^{2m}(\la(1-i/N)) \times (L_i \cup
P_{i+1})\right)\;\;\cup\;\;\dots\vspace{7mm}\\ \hspace{13mm}\cup\;\;
\left(B^{2m}(\la/N)\times P_N\right).
\end{eqnarray*}
Thus the interval $L_1$ for example is longer than we need: $G$ maps it into a
hypersurface which disjoins the whole ball $B^{2m}(\la)$ while we only need the
hypersurface to disjoin  the subball $B^{2m}(\la(1-1/N))$ from $B^{2m}(\la)$. Because
the hypersurface comes from the regular isotopy $\si_s$, we only need the part of it
corresponding to the isotopy $\si_s, 0\leq s\leq 1-1/N$, which sits over an interval of
length $1 - 1/N$.
 
Here are the details.  Let us change the proportions in
the set $Y_N$ by moving the points $ b_i$, keeping
the $a_i$ fixed.  Thus, for $0 \leq s \leq 1$, put $b_i^s = b_i + s$, and let
$$
Y_N^s =  P_1^s \cup L_1^s \cup\dots\cup P_N^s,
$$
where the rectangles $P_i^s$ have area $ (1/N +s)\la$ and the intervals $L_i$ have length $1-s$.
Then, for each $s$, $B^{2m+2}(\la)$ embeds in a neighborhood of
\begin{eqnarray*}
W^s = \left(B^{2m}(\la)\times P_1^s\right) \;\;\cup\;\;\left( B^{2m}(\la(1-1/N-s)) \times
(L_1^s \cup P_2^s)\right) \;\;\cup\;\;\dots  \vspace{7mm}\\
\hspace{10mm} \cup \;\;\left(B^{2m}(\la(1-i/N -is)) \times (L_i^s \cup
P_{i+1}^s)\right)\;\;\cup\;\;\dots\vspace{7mm}\\ \hspace{13mm}\cup\;\;
\left(B^{2m}(\la(1/N - (N-1)s))\times P_N^s\right).
\end{eqnarray*}
Here, the convention is that $B^{2m}(\ka)$ denotes a point when $\ka= 0$ and  is empty
for $\ka \leq 0$. Note that when
$s = s_0 = 1-1/N$, $W^{s_0}$ is just $B^{2m}(\la)\times P_1^{s_0} $ embedded by a
standard  embedding  and the isotopy may be ended.
The length of $L_i^{s_0}$ is then $\frac 1N > \de$:  the lengths of $L_i^s$
are thus larger than $\de$ during the whole isotopy.
 
We now map $W^s$ to $V\times S$ by the obvious map $G^s$, which, for each 
$0 \leq s \leq 1 - 1/N$ is constructed from the  isotopy $\si_u, \, 0\leq u \leq 1-s$.  
It is not hard to check that the
$\de$\/-regularity of $\si_s$ implies that
$\im\, G^s$ is disjoined by the translation $\tau_S$ along $S$.
Since we may choose everything to
vary smoothly with $s$, the result follows.  \QED
 
\noindent
{\bf Proof of Lemma~\ref{unwrap}}$\quad$  We  construct a family of balls $B_t'$
and isotopies $\si^t_s$  from the unwrapped ball $B'(c) = B(B,\phi)$ to a
ball and isotopy which are the lifts of a ball and strictly disjoining isotopy in $\R^2$.   The
first part of the isotopy is simple: we just decrease the capacity
of the unwrapped  ball  $B'$ from $c$ to $\la$ by shrinking the initial ball $B(c)$ while keeping
the initial disjoining diffeotopy $\phi_s$ fixed.
Note that $\tau_s$ does not change through this isotopy.
 
      Consider now the ball $B'(\la)$ obtained at the end of this isotopy.
It is clear by the definition of the unwrapped construction that the
following diagram commutes:
$$
\begin{array}{ccccc}
B^{2n+4}(\la) & \stackrel{i}{\hookrightarrow} & B^{2n+2}(\la) \times \Nn(Y_N(\la)) &
\stackrel{\tilde G_g}{\hookrightarrow} & M' \times \R^2  \vspace{7mm}\\
              &                   & \downarrow (\pi \circ g) \times id
&                     & \downarrow (\pi \times id) \vspace{7mm}\\
    &
 & B^{2}(\la) \times \Nn(Y_N(\la))    & \stackrel{G_j}{\hookrightarrow}
& \R^2 \times \R^2 \end{array}
$$
where $g$ is the initial embedding $B^{2n}(\la) \hookrightarrow B'(\la)$, $j$ is
 the inclusion
$B^2(\la)\hookrightarrow \R^2$, $\pi$ is the projection onto
$\R^2$,
and we write $\tilde G_g$, $G_j$ for the corresponding maps which
give the unwrapped balls. It follows that:
\medskip
 
\noindent
(i) $\quad$  any
isotopy $ G_t, 0 \leq t \leq t_0$
beginning with $ G_0 = G_j$
lifts to an isotopy $\tilde{G}_t$ beginning with $\tilde{G}_0 =
\tilde G_g$.
 
\medskip
\noindent
(ii) $\quad$  If, for each $t$, $ \rho_s^t$ is a $1$\/-periodic
diffeotopy of $\R^4$ which strictly
disjoins $\im\, ( G_t)$, then
its pull-back $\tilde{\rho}_s^t$ to $M \times \R^4$ also
strictly disjoins $\im\,(\tilde{G}_t)$
and has the same energy.
 
\medskip
\noindent
(iii) $\quad$  Finally, if at time $t_0$, $G_{t_0}$ is a product
$$
f_1 \times f_2:B^{2}(\la) \times N(Y_N(\la))  \hookrightarrow \R^2 \times \R^2
$$
and $\rho_s^{t_0}$ is (the pull-back of) a translation $\psi_s$ in
the second $\R^2$\/-factor such
that the pair
$(f_2, \psi_s)$ is $\de$\/-regular, then $\tilde{G}_{t_0}$ is again a
product whose second factor is
$\de$\/-regular under $\psi_s$.  It follows that the pair
$(\tilde{G}_{t_0} \circ i, \tilde{\rho}_s^{t_0})$
is also $\de$\/-regular by Lemma~\ref{regproduct}.
\medskip
 
    This reduces the proof to the construction of the $4$\/-dimensional isotopy $G_t$.
This may be easily defined by using the fact that the translation $\phi_s$ which
strictly disjoins $B^2(\la)$ in $\R^2$ is generated by an  {\it autonomous} Hamiltonian
$H$ (of total variation $\la + \eps$): for the moment, let us forget the
(technical) requirement that our isotopies should be constant when the time $s$
is near
an integer value, and let us suppose that the isotopy $\si_s$ is a $1$\/-parameter
 group.
Then, for each $0 \leq t \leq 1$, set
$$
\si_s^{t} = \si_{(1-t)s},
$$
and consider the corresponding unwrapped ball $B(B^2(\la), \si^{t}) = \im\,
G_{\si^{t}}$.  Because the energy $e(\si^{t})$ of $\si^{t}$ is $(1-t)e(\si) =
(1-t)(\la +
\eps)$, the ball  $B(B^2(\la), \si^{t})$ sits over a strip $S^t=\coprod_i
P_i\cup\coprod_iQ_i^t$  where the rectanges $Q_i^t$ have area $(1-t)(\la + \eps)
$.
Therefore the translation  $\tau^t$ which moves $S^t$ through the distance $1 +1/N$
may be extended as in Lemma~\ref{extend} to an isotopy $\tau^t_s$ of energy
$$
\frac{\la}{N} +
e(\si^{t}) + \eps/2 = \frac{\la}{N} + (1-t)(\la + \eps) + \eps/2.
$$
Note that $\tau^t$  does not disjoin $B(B^2(\la), \si^{t})$ from itself, since
 $\si^{t}$
does not disjoin $B^2(\la)$ at time $s=1$.  However, if we follow  $\tau^t$ (which is a
movement in the second $\R^2$ factor) with the translation $\si_1^{t+1} = \si_{-t}$ in the
first $\R^2$ factor, we do get a disjoining isotopy.  Thus, the isotopy
$$
\rho^t_s = \si_{-ts}\circ \tau^{t}_s
$$
disjoins $B(B^2(\la), \si^{t})$ from itself at time $s=1$, and may be extended
 to an
isotopy of $\R^4$ with total energy
\begin{eqnarray*}
\frac{\la}{N} +
e(\si^{t}) + e(\si^{t+1}) + \eps &=& \frac{\la}{N} + (1-t)(\la + \eps) + t(\la
+ \eps) +\eps \\
&=&  \frac{\la}{N} + \la + 2\eps.
\end{eqnarray*}
  Taking $G^t =
G_{\si^{t}}$, we can therefore satisfy (i) and (ii) above.  Observe that when
$t=1$,
$G_1$ is simply the inclusion of $B^2(\la) \times N(Y_N(\la))$ and 
$\rho_s^1$ is (isotopic to) the
disjoining translation $(\si_s)^{-1}$ in the first $\R^2$-direction. It is easy to see that there 
is a family of symplectomorphisms $\be_t, 1\leq t \leq 2$ of $\R^4$ which begins with the identity
such that, for each $t
\in [1,2]$, the conjugate isotopies $\be_t^{-1} \circ\rho_s^1\circ \be_t$ strictly disjoin
 $\im\,(G_1)$  and such that  the final isotopy $\psi_s = \be_2^{-1}
\circ\rho_s^1\circ \be_2$ has the form  of a translation in the $u$\/-direction
$$
(x,y, u,v) \mapsto (x,y,u, v + \alpha_s(u))
$$
that can  be chosen  $\de$\/-regular on the second factor of $G_1$.
The pair $(G_1, \psi_s)$ then
satisfies the conditions in (iii) above. 

It remains to take into account the requirement that
our isotopies should be constant when the time $s$ is near an integral value.
But it is  sufficient to choose the initial translation $\phi_s$, that disjoins $B^2(\la)$
in $\R^2$, generated by an Hamiltonian of the form $H_s = hH$ where $H$ is autonomous and 
$h:\R \to \R$ is a bump function
equal to $1$ everywhere except on a $\mu$\/-neighborhood of each
integer: choosing $\mu$ small enough with respect to $\de$, we may clearly arrange that the above
argument holds.    
\QED


\begin{thebibliography}{99999}
 
\bibitem{EH1}  Y. Eliashberg and H. Hofer:  {\it An energy-capacity inequality for the
symplectic holonomy of hypersurfaces flat at infinity}, Preprint 1992.
 
\bibitem{EH2}  Y. Eliashberg and H. Hofer:  {\it Towards the definition of symplectic
boundary}, Preprint 1993.
 
 
\bibitem{EP}  Y. Eliashberg and L. Polterovich:  {\it  Biinvariant metrics on the group
of Hamiltonian diffeomorphisms}, preprint 1991
 
\bibitem{FHW} A. Floer, H. Hofer and K. Wysocki: {\it Applications of Symplectic
Homology I}, preprint 1992.

\bibitem{GR} Gromov, M.:  {\it Pseudo-holomorphic curves in symplectic manifolds},
Invent. Math. {\bf 82},  1985, 307-347.
 
 
\bibitem{H1}  H. Hofer:  {\it  On the topological properties of symplectic maps}, Proc.
Royal Soc. Edinburgh {\bf 115} A (1990), 25-38.
 
\bibitem{H2}  H. Hofer:  {\it Estimates for the energy of a symplectic map},Comment.
Math. Helv. {\bf 68} (1993), 48-72.
 
\bibitem{H3}  H. Hofer:  {\it Symplectic Capacities}, Durham Conference, ed:
Donaldson and Thomas, London Math Soc, 1992.
 
\bibitem{LA}  F. Lalonde, {\it Isotopy of symplectic balls, Gromov's radius, and the
structure of ruled symplectic 4-manifolds}, preprint 1992.
 
\bibitem{ELL} D. McDuff: {\it Elliptic methods in symplectic topology}, Bull. Amer.
Math Soc. {\bf 23} (1990), 311-358.
 
\bibitem{LB} D. McDuff: {\it Singularities of $J$-holomorphic curves}, Journ. Geom.
Anal, 1992
 
 \bibitem{CTB} D. McDuff:  {\it Symplectic manifolds with boundaries of
contact-type boundaries }, Invent. Math. {\bf 103} (1991), 651-671.
 
\bibitem{UNIQ} D. McDuff: {\it Remarks on the uniqueness of blowing up}, to appear in Proceedings of Conference in Warwick 1990, ed Salamon.
 
\bibitem{MP} D. McDuff and L. Polterovich, {\it Symplectic packings and Algebraic
Geometry}, preprint 1992

\bibitem{MS} D. McDuff and D. Salamon, {\it Notes on $J$-holomorphic curves },
 Stony Brook preprint 1993.
 
\bibitem{POL}  L. Polterovich: {\it Symplectic displacement energy for Lagrangian
submanifolds}, preprint 1991.
 
\bibitem{TR} L. Traynor, in preparation.
 
\bibitem{V}  C. Viterbo:  {\it Symplectic topology as the geometry of generating
functions}, Math. Ann {\bf 292} (1992), 685-710.
 
\end{thebibliography}
\end{document}